\documentclass[12pt, reqno]{amsart}
\usepackage{amsmath, amsthm, amscd, amsfonts, amssymb, graphicx, color}
\usepackage[bookmarksnumbered, colorlinks=false, plainpages]{hyperref}

\newtheorem{theorem}{Theorem}[subsection]
\newtheorem{lemma}[theorem]{Lemma}
\newtheorem{proposition}[theorem]{Proposition}
\newtheorem{corollary}[theorem]{Corollary}
\theoremstyle{definition}
\newtheorem{definition}[theorem]{Definition}
\newtheorem{example}[theorem]{Example}

\theoremstyle{remark}
\newtheorem{remark}[theorem]{Remark}
\numberwithin{equation}{section}

\allowdisplaybreaks
\begin{document}

\title[$E$-Bessel sequences and $E$-multipliers in Hilbert spaces]
{$E$-Bessel sequences and $E$-multipliers in Hilbert spaces}

\author[H.Hedayatirad]{Hassan Hedayatirad}
\address{Hassan Hedayatirad \\ Department of Mathematics and Computer
Sciences, Hakim Sabzevari University, Sabzevar, P.O. Box 397, IRAN}
\email{ \rm hasan.hedayatirad@hsu.ac.ir; hassan.hedayatirad67@gmail.com}
\author[T.L. Shateri]{Tayebe Lal Shateri }
\address{Tayebe Lal Shateri \\ Department of Mathematics and Computer
Sciences, Hakim Sabzevari University, Sabzevar, P.O. Box 397, IRAN}
\email{ \rm  t.shateri@hsu.ac.ir; t.shateri@gmail.com}
\thanks{*The corresponding author:
t.shateri@hsu.ac.ir; t.shateri@gmail.com (Tayebe Lal Shateri)}
 \subjclass[2010] {Primary 42C15;
Secondary 54D55.} \keywords{$E$-frame, Hilbert space, Multiplier, $E$-multiplier.}
 \maketitle
\begin{abstract}
$E$-frames are a new generalization for the concept of frames for $\mathcal{H}$, where $E$ is an infinite invertible complex matrix mapping on $\bigoplus_{n=1}^{\infty}\mathcal{H}$. This article is dedicated to investigating some notions related to $E$-Bessel sequences and $E$-multipliers. A Multiplier is an operator created by frame-like analysis, multiplication with a fixed sequence, called the symbol, and synthesis. In this article, we introduce the notion of $E$-multipliers, which generalizes multipliers for $E$-sequences and study their properties, including boundedness and invertibility.

\end{abstract}
\section{Introduction}
Suppose that $\mathcal{H}$ is a separable Hilbert space. Recall that a sequence $\left\lbrace\psi_k\right\rbrace_{k=1}^{\infty}$ is a frame for $\mathcal{H}$ if and only if there exists $0<A_\Psi\leq B_\Psi<\infty$ such that for each $f\in\mathcal{H}$ 
\begin{equation}\label{frame}
A_\Psi\left\Vert f\right\Vert^2\leq\sum_{k=1}^{\infty}\left\vert\left\langle f,\psi_k\right\rangle\right\vert^2\leq B_\Psi\left\Vert f\right\Vert^2.
\end{equation}
Consider an infinite invertible matrix mapping $E$ on $\bigoplus_{n=1}^{\infty}\mathcal{H}=\left\lbrace\left\lbrace f_n\right\rbrace_{n=1}^{\infty}\;,\;\sum_{n=1}^{\infty}\left\Vert f_n\right\Vert^2<\infty\right\rbrace$. A sequence $\left\lbrace\psi_k\right\rbrace_{k=1}^{\infty}$ is called an $E$-frame if and only if $E\left\lbrace \psi_k\right\rbrace_{k=1}^{\infty}$ is a well defined sequence in $\mathcal{H}$ and some positive numbers $A\leq B$ exist such that
\begin{equation}\label{eframe}
A\left\Vert f\right\Vert^2\leq\sum_{n=1}^{\infty}\left\vert\left\langle f,\left(E\left\lbrace \psi_k\right\rbrace_{k=1}^{\infty}\right)_n\right\rangle\right\vert^2\leq B\left\Vert f\right\Vert^2,
\end{equation}
for all $f\in\mathcal{H}$ \cite{TAL}. If only the right inequality holds in \eqref{frame} and \eqref{eframe}, we say $\left\lbrace\psi_k\right\rbrace_{k=1}^{\infty}$ is Bessel and $E$-Bessel respectively. Associated to each Bessel sequence $\Psi=\left\lbrace\psi_k\right\rbrace_{k=1}^{\infty}$, the synthesis (pre-frame) operator $T_\Psi:\ell^2(\mathbb{N})\longrightarrow\mathcal{H}\;;\;T_\Psi\left\lbrace c_k\right\rbrace_{k=1}^{\infty}=\sum_{k=1}^{\infty}c_k\psi_k$ is well defined and bounded. The adjoint operator of $T$, is called the analysis operator and defined by $T^\ast_\Psi:\mathcal{H}\longrightarrow\ell^2(\mathbb{N})\;;\;T^\ast_\Psi f=\left\lbrace\left\langle f,\psi_k\right\rangle\right\rbrace_{k=1}^{\infty}.$ The operator $S=T_\Psi T^\ast_\Psi$ is the frame operator which is self adjoint, positive and invertible \cite{ch}.

Similarly, $\Psi=\left\lbrace\psi_k\right\rbrace_{k=1}^{\infty}$ is an $E$-Bessel sequence if and only if the synthesis (pre-$E$-frame) operator $T_{E}:\ell^2(\mathbb{N})\longrightarrow\mathcal{H}\;;\;T_E\left\lbrace c_n\right\rbrace_{n=1}^{\infty}=\sum_{n=1}^{\infty}c_n\left(E\left\lbrace\psi_k\right\rbrace_{k=1}^{\infty}\right)_n$ is well defined and bounded. The adjoint of $T_E$ is defined by $T^\ast_E:\mathcal{H}\longrightarrow\ell^2(\mathbb{N})\;;\;T^\ast_E f=\left\lbrace\left\langle f,\left(E\left\lbrace\psi_k\right\rbrace_{k=1}^{\infty}\right)_n\right\rangle\right\rbrace_{n=1}^{\infty}.$ Finally, $S_E=T_ET^\ast_E$ is a self adjoint, positive and invertible operator which is called the $E$-frame operator \cite{TAL}.

Balazs introduced the Bessel and frame multipliers for Hilbert spaces \cite{BAL1}. There are numerous applications of this kind of operators.  Such operators find application in psychoacoustics \cite{BAL3}, virtual acoustics \cite{MAJ1}, denoising \cite{MAJ2}.\\
For any sequence $\Psi=\left\lbrace\psi_k\right\rbrace_{k=1}^{\infty}$ and $\Phi=\left\lbrace\phi_k\right\rbrace_{k=1}^{\infty}$ in $\mathcal{H}$ and any sequence of complex numbers $m=\left\lbrace m_k\right\rbrace_{k=1}^{\infty}$ (called symbol), the operator $M_{m,\Psi,\Phi}$, given by
\begin{equation}\label{multi}
M_{m,\Psi,\Phi}f=\sum_{k=1}^{\infty}m_k\left\langle f,\phi_k\right\rangle\psi_k,
\end{equation}
is called a multiplier \cite{BAL1}. We say $M_{m,\Psi,\Phi}$ is well defined on $\mathcal{H}$ if the series in \eqref{multi} converges for all $f\in\mathcal{H}$. If $\sum_{k=1}^{\infty}m_k\left\langle f,\phi_k\right\rangle\psi_k$ is unconditionally convergent, then we say $M_{m,\Psi,\Phi}$ is unconditionally convergent. Depending on $\Psi$, $\Phi$ and $m$, the corresponding multiplier might not be well defined, it might be well defined but not unconditionally convergent and it might be unconditionally convergent \cite{BAL2}.

In this paper, we first show that any infinite matrix mapping on $\bigoplus_{n=1}^{\infty}\mathcal{H}$ can be considered as an operator on $\ell^2(\mathbb{N})$, under certain conditions. Using this, we reintroduce $T_E$ and $T^\ast_E$ in terms of $T_\Psi$ and $T^\ast_\Psi$ where $\Psi$ is a fixed Bessel sequence in $\mathcal{H}$. In section 1, we try to reintroduce the operators associated with $E$-frames in terms of the frame's analysis and synthesis operators, under certain conditions. In section 2, we introduce the notion of $E$-multipliers which is the extension of the concept of multipliers for $E$-sequences. We review the behavior of these operators when the parameters are changing. Also, we investigate, when an $E$-multiplier is bounded or invertible.

As mentioned above, in \cite{TAL}, it is assumed that $E$ is an infinite matrix that defines a mapping on $\bigoplus_{n=1}^{\infty}\mathcal{H}$. In our definition of the $E$-frame, we remove this condition and will only include it as an assumption if necessary.
\section{Main Results}
Throughout this section, $\mathcal{H}$ is a separable Hilbert space, $\left\lbrace e_n\right\rbrace_{n=1}^{\infty}$ is an orthonormal basis and $E$ is an invertible infinite complex matrix. First, we want to check some conditions to transmute a matrix mapping on $\bigoplus_{n=1}^{\infty}\mathcal{H}$ into a matrix mapping on $\ell^2(\mathbb{N})$.
\subsection{$E$-Bessel sequences}
\begin{remark}\label{rem1}
Let $E$ be an infinite complex matrix mapping on $\bigoplus_{n=1}^{\infty}\mathcal{H}$ defined by
\begin{equation*}
E:\bigoplus_{n=1}^{\infty}\mathcal{H}\longrightarrow\bigoplus_{n=1}^{\infty}\mathcal{H}\;;\;E\lbrace f_k\rbrace_{k=1}^{\infty}=\left\lbrace\sum_{k=1}^{\infty}E_{n,k}f_k\right\rbrace_{n=1}^{\infty},
\end{equation*}
which satisfies 
\begin{equation}\label{eq1}
\sum_{n=1}^{\infty}\sum_{k=1}^{\infty}\left\vert E_{n,k}\right\vert^2<\infty.
\end{equation}
Then $E$ is bounded. In fact,
\begin{align}\label{Ebdd}
\Vert E\Vert^2=\sup_{\sum_k\left\Vert f_k\right\Vert^2\leq1}\left\Vert E\left\lbrace f_k\right\rbrace_{k=1}^{\infty}\right\Vert_{\bigoplus\mathcal{H}}^2&=\sup_{\sum_k\left\Vert f_k\right\Vert^2\leq1}\sum_{n=1}^{\infty}\left\Vert\sum_{k=1}^{\infty}E_{n,k}f_k\right\Vert^2\\&\leq\sup_{\sum_k\left\Vert f_k\right\Vert^2\leq1}\sum_{n=1}^{\infty}\left(\sum_{k=1}^{\infty}\left\vert E_{n,k}\right\vert\left\Vert f_k\right\Vert\right)^2\nonumber\\&\leq\sup_{\sum_k\left\Vert f_k\right\Vert^2\leq1}\sum_{n=1}^{\infty}\sum_{k=1}^{\infty}\left\vert E_{n,k}\right\vert^2\sum_{k=1}^{\infty}\left\Vert f_k\right\Vert^2<\infty\nonumber.
\end{align}
Now, we use this to directly prove that $\lbrace E_{n,k}\rbrace_{k=1}^{\infty}\in\ell^2(\mathbb{N})$ for all $n\in\mathbb{N}$. Consider an orthonormal basis $\lbrace e_k\rbrace_{k=1}^{\infty}$ for $\mathcal{H}$ and for a fixed $j\in\mathbb{N}$, suppose that $\lbrace \delta_{j,k}\rbrace_{k=1}^{\infty}$ is a sequence in $\mathcal{H}$ defined by
\begin{equation*}
\delta_{j,k}=\begin{cases} 
e_j&k=j\\0&k\neq j.
\end{cases}
\end{equation*}
It is clear that $\lbrace \delta_{j,k}\rbrace_{k=1}^{\infty}$ belongs to $\bigoplus_{n=1}^{\infty}\mathcal{H}$. Furthermore,
\begin{align*}
E\left\lbrace\delta_{j,k}\right\rbrace_{k=1}^{\infty}=\left\lbrace\sum_{k=1}^{\infty}E_{n,k}\delta_{j,k}\right\rbrace_{n=1}^{\infty}=\left\lbrace E_{n,j}e_j\right\rbrace_{n=1}^{\infty}.
\end{align*}
Therefore, \eqref{Ebdd} implies that
\begin{align}\label{Ecol}
\sum_{n=1}^{\infty}\left\vert E_{n,j}\right\vert^2=\sum_{n=1}^{\infty}\left\Vert E_{n,j}e_j\right\Vert^2=\left\Vert\left\lbrace E_{n,j}e_j\right\rbrace_{n=1}^{\infty}\right\Vert_{\bigoplus\mathcal{H}}^2=\left\Vert E\left\lbrace\delta_{j,k} \right\rbrace_{k=1}^{\infty}\right\Vert_{\bigoplus\mathcal{H}}^2\leq\Vert E\Vert.
\end{align}
(\ref{Ecol}) shows that each columns of $E$ belongs to $\ell^2(\mathbb{N})$. An analogous argument applying on $E^{\ast}$ proves that $\left\lbrace E_{n,k}\right\rbrace_{k=1}^{\infty}\in\ell^2(\mathbb{N})$.

Now we are going to define an operator on $\ell^2(\mathbb{N})$ using $E$. We will denote it again by $E$ and define it as follows
\begin{equation*}
E:\ell^2(\mathbb{N})\longrightarrow\ell^2(\mathbb{N})\;;\;E\lbrace c_k\rbrace_{k=1}^{\infty}=\left\lbrace\sum_{k=1}^{\infty}E_{n,k}c_k\right\rbrace_{n=1}^{\infty}.
\end{equation*}
The above discussion together with \eqref{eq1}, makes $E$ well defined. Using reasoning similar to that in \eqref{Ebdd}, we can conclude that $E$ is bounded.
Also, one can easily shows that 
\begin{align*}
E^t\lbrace c_k\rbrace_{k=1}^{\infty}=\left\lbrace\sum_{k=1}^{\infty}E_{k,n}c_k\right\rbrace_{n=1}^{\infty}\qquad\text{and}\qquad \overline{E}\lbrace c_k\rbrace_{k=1}^{\infty}=\left\lbrace\sum_{k=1}^{\infty}\overline{E_{n,k}}c_k\right\rbrace_{n=1}^{\infty},
\end{align*} 
are well defined matrix operators on $\ell^2(\mathbb{N})$.
\end{remark}
\begin{proposition}
Suppose that $E$ is an infinite complex matrix mapping on $\bigoplus_{n=1}^{\infty}\mathcal{H}$ satisfies \eqref{eq1}, $\Psi=\lbrace\psi_k\rbrace_{k=1}^{\infty}$ is an $E$-Bessel sequence on $\mathcal{H}$ and $m=\lbrace m_k\rbrace_{k=1}^{\infty}$ is a scalar sequence such that $m\Psi$ is Bessel. Then $m\Psi$ is an $E$-Bessel sequence for $\mathcal{H}$.
\end{proposition}
\begin{proof}
We use Remark \ref{rem1} to treat $E$ as an operator on $\bigoplus_{n=1}^{\infty}\mathcal{H}$ or $\ell^2(\mathbb{N})$, depending on the context. First note that $E$ is bounded and the rows of $E$ belongs to $\ell^2(\mathbb{N})$ by Remark \ref{rem1}. Hence $\sum_{k=1}^{\infty}E_{n,k}m_k\psi_k$ is well defined for all $n\in\mathbb{N}$. Moreover for given $f\in\mathcal{H}$
\begin{align*}
\sum_{n=1}^{\infty}\left\vert\left\langle f,\left(E\left\lbrace m_k\psi_k\right\rbrace_{k=1}^{\infty}\right)_n\right\rangle\right\vert^2&=\sum_{n=1}^{\infty}\left\vert\left\langle f,\sum_{k=1}^{\infty}E_{n,k}m_k\psi_k\right\rangle\right\vert^2\\&=\sum_{n=1}^{\infty}\left\vert\sum_{k=1}^{\infty}E_{n,k}\left\langle f,m_k\psi_k\right\rangle\right\vert^2\\&=\sum_{n=1}^{\infty}\left\vert\left(E\left\lbrace\left\langle f,m_k\psi_k\right\rangle\right\rbrace_{k=1}^{\infty}\right)_n\right\vert^2\\&=\left\Vert E\left\lbrace\left\langle f,m_k\psi_k\right\rangle\right\rbrace_{k=1}^{\infty}\right\Vert_{\ell^2}^2\leq\left\Vert E\right\Vert^2B\left\Vert f\right\Vert^2,
\end{align*}
where $B$ is the Bessel bound for $m\Psi$.
\end{proof}
\begin{corollary}
Suppose that $E$ is an infinite complex matrix mapping on $\bigoplus_{n=1}^{\infty}\mathcal{H}$ satisfies \eqref{eq1}, $\Psi=\lbrace\psi_k\rbrace_{k=1}^{\infty}$ is an $E$-Bessel sequence on $\mathcal{H}$ and $m=\lbrace m_k\rbrace_{k=1}^{\infty}$ belongs to $\ell^\infty$. Then $m\Psi$ is an $E$-Bessel sequence for $\mathcal{H}$.
\end{corollary}
As preparation for Theorem \ref{th2.1}, we show:
\begin{proposition}\label{pro1}
Let $E$ be an infinite complex matrix mapping on $\bigoplus_{n=1}^{\infty}\mathcal{H}$ which satisfies \eqref{eq1}. Suppose that $\lbrace f_k\rbrace_{k=1}^{\infty}\in\bigoplus_{n=1}^{\infty}\mathcal{H}$ is an $E$-Bessel sequence. Then for each $\lbrace c_n\rbrace_{n=1}^{\infty}\in\ell^2(\mathbb{N})$, the series 
$$\sum_{n=1}^{\infty}\sum_{k=1}^{\infty} c_nE_{n,k}f_k$$ is absolutely convergent.
\end{proposition}
\begin{proof}
\begin{align*}
\sum_{n=1}^{\infty}\sum_{k=1}^{\infty}\left\Vert c_nE_{n,k}f_k\right\Vert&=\sum_{n=1}^{\infty}\left\vert c_n\right\vert\sum_{k=1}^{\infty}\left\vert E_{n,k}\right\vert\left\Vert f_k\right\Vert\\&\leq\left\lbrace\sum_{n=1}^{\infty}\left\vert c_n\right\vert^2\right\rbrace^{\frac{1}{2}}\left\lbrace\sum_{n=1}^{\infty}\left(\sum_{k=1}^{\infty}\left\vert E_{n,k}\right\vert\left\Vert f_k\right\Vert\right)^2\right\rbrace^{\frac{1}{2}}\\&\leq\left\lbrace\sum_{n=1}^{\infty}\left\vert c_n\right\vert^2\right\rbrace^{\frac{1}{2}}\left\lbrace\sum_{n=1}^{\infty}\sum_{k=1}^{\infty}\left\vert E_{n,k}\right\vert^2\sum_{k=1}^{\infty}\left\Vert f_k\right\Vert^2\right\rbrace^{\frac{1}{2}}\\&=\left\lbrace\sum_{n=1}^{\infty}\left\vert c_n\right\vert^2\right\rbrace^{\frac{1}{2}}\left\lbrace\sum_{k=1}^{\infty}\left\Vert f_k\right\Vert^2\right\rbrace^{\frac{1}{2}}\left\lbrace\sum_{n=1}^{\infty}\sum_{k=1}^{\infty}\left\vert E_{n,k}\right\vert^2\right\rbrace^{\frac{1}{2}}<\infty.
\end{align*}
\end{proof}
\begin{theorem}\label{th2.1}
Let $E$ be an infinite complex matrix mapping on $\bigoplus_{n=1}^{\infty}\mathcal{H}$ which satisfies \eqref{eq1}. Suppose that $\lbrace f_k\rbrace_{k=1}^{\infty}\in\bigoplus_{n=1}^{\infty}\mathcal{H}$ is an $E$-Bessel sequence. Then 
\begin{itemize}
\item[(i)] For each sequence $\lbrace c_n\rbrace_{n=1}^{\infty}\in\ell^2(\mathbb{N})$, 
\begin{equation*}
T_{E}\left\lbrace c_n\right\rbrace_{n=1}^{\infty}=T\left(E^t\left\lbrace c_n\right\rbrace_{n=1}^{\infty}\right).
\end{equation*}
\item[(ii)] For each $f\in\mathcal{H}$, 
\begin{equation*}
T_{E}^{\ast}f=\overline{E}T^{\ast}f.
\end{equation*}
\end{itemize}
\end{theorem}
\begin{proof}
Using Remark \ref{rem1} and Proposition \ref{pro1} and applying the Fubini's theorem we have
\begin{align}
T_{E}\left\lbrace c_n\right\rbrace_{n=1}^{\infty}&=\sum_{n=1}^{\infty}c_n\left(E\left\lbrace f_k\right\rbrace_{k=1}^{\infty}\right)_n\\&=\sum_{n=1}^{\infty}\sum_{k=1}^{\infty
}c_{n}E_{n,k}f_k\nonumber\\&=\sum_{k=1}^{\infty
}\left(\sum_{n=1}^{\infty}c_{n}E_{n,k}\right)f_k=T\left\lbrace\sum_{n=1}^{\infty}c_{n}E_{n,k}\right\rbrace_{k=1}^{\infty}=T\left(E^{t}\left\lbrace c_n\right\rbrace_{n=1}^{\infty}\right),\nonumber
\end{align}
where $E^{t}$ is the transpose of $E$.
\begin{align*}
T_{E}^{\ast}f=\left\lbrace\left\langle f,\left(E\left\lbrace f_k\right\rbrace_{k=1}^{\infty}\right)_n\right\rangle\right\rbrace_{n=1}^{\infty}&=\left\lbrace\left\langle f,\sum_{k=1}^{\infty}E_{n,k}f_k\right\rangle\right\rbrace_{n=1}^{\infty}\\&=\left\lbrace\sum_{k=1}^{\infty}\overline{E_{n,k}}\left\langle f,f_k\right\rangle\right\rbrace_{n=1}^{\infty}\\&=\overline{E}\left\lbrace\left\langle f,f_k\right\rangle\right\rbrace_{k=1}^{\infty}=\overline{E}T^{\ast}f.
\end{align*}
Note that $\lbrace f_k\rbrace_{k=1}^{\infty}$ is a Bessel sequence and so $\left\lbrace\left\langle f,f_k\right\rangle\right\rbrace_{k=1}^{\infty}$ belongs to $\ell^2(\mathbb{N})$.
\end{proof}
\subsection{E-multiplier}
In the sequel, the concept of $E$-multiplier operators for $E$-Bessel sequences will be introduced and some of their properties will be shown. These operators are defined by a
fixed multiplication pattern which is inserted between the
analysis and synthesis operators.
\begin{definition}
Let $E$ be an infinite complex matrix. We say a sequence $\lbrace f_k\rbrace_{k=1}^{\infty}$ is an $E$-sequence in $\mathcal{H}$ if, its $E$-transform, i.e. the sequence 
\begin{equation*}
E\lbrace f_k\rbrace_{k=1}^{\infty}=\left\lbrace\sum_{k=1}^{\infty}E_{n,k}f_k\right\rbrace_{n=1}^{\infty},
\end{equation*}
is a well defined sequence in $\mathcal{H}$.
\end{definition}
The complex sequence $m=\lbrace m_n\rbrace_{n=1}^{\infty}$ is called semi-normalized if $0<\inf_n\vert m_n\vert\leq\sup_n\vert m_n\vert<\infty$. The sequence $\Phi=\lbrace\phi_k\rbrace_{k=1}^{\infty}$ is called norm-bounded below (resp. norm-bounded above) if $\inf_k\Vert\phi_k\Vert>0$ (resp. $\sup_k\Vert\phi_k\Vert<\infty$) and $\Phi$ is called semi-normalized if $0<\inf_k\Vert\phi_k\Vert\leq\sup_k\Vert\phi_k\Vert<\infty$.
\begin{definition}
Consider two separable Hilbert spaces $\mathcal{H}_1$ and $\mathcal{H}_2$ and suppose that $E_1$ and $E_2$ are infinite complex matrices. For any $E_1$-sequence $\Phi=\lbrace\phi_k\rbrace_{k=1}^{\infty}$ in $\mathcal{H}_1$ and $E_2$-sequence $\Psi=\lbrace\psi_k\rbrace_{k=1}^{\infty}$ in $\mathcal{H}_2$ and any sequence of complex number $m=\lbrace m_n\rbrace_{n=1}^{\infty}$ (we call it the symbol), the operator $M_{m,\Psi,\Phi}^{\left(E_1,E_2\right)}$ from $\mathcal{H}_1$ to $\mathcal{H}_2$ defined by
\begin{equation}\label{eqmlt}
M_{m,\Psi,\Phi}^{\left(E_1,E_2\right)}f=\sum_{n=1}^{\infty}m_n\left\langle f,\left(E_1\left\lbrace\phi_k\right\rbrace_{k=1}^{\infty}\right)_n\right\rangle\left(E_2\left\lbrace\psi_k\right\rbrace_{k=1}^{\infty}\right)_n,
\end{equation}
is a $(E_1,E_2)$-Multiplier.
\end{definition}
Depending on $m$, $\Psi$, $\phi$, $E_1$ and $E_2$, the series in \eqref{eqmlt} might not converges for some $f\in\mathcal{H}_1$. 
\begin{remark}\label{rem2}
As in \cite{BAL1}, we use the notation $\mathcal{M}_m$ for the operator $\mathcal{M}_m:\ell^2(\mathbb{N})\longrightarrow\ell^2(\mathbb{N})$ with $m=\{m_k\}\in\ell^p (p>0)$ defined by $\mathcal{M}_m\lbrace c_k\rbrace_{k=1}^\infty=\lbrace m_kc_k\rbrace_{k=1}^\infty$. If $\Phi=\lbrace\phi_k\rbrace_{k=1}^{\infty}$ and $\Psi=\lbrace\psi_k\rbrace_{k=1}^{\infty}$ are $E_1$-Bessel and $E_2$-Bessel sequences for Hilbert spaces $\mathcal{H}_1$ and $\mathcal{H}_2$  with bounds $B$ and $B^{'}$ respectively, then we can write 
\begin{equation}\label{eqmlt2}
M_{m,\Psi,\Phi}^{\left(E_1,E_2\right)}=T_{E\Psi}\mathcal{M}_mT^\ast_{E\Phi}.
\end{equation}
\end{remark}
We will show below that by considering the conditions of Remark \ref{rem2}, if $m\in\ell^\infty$, then $M_{m,\Psi,\Phi}^{\left(E_1,E_2\right)}$ is well defined and bounded. 
\begin{theorem}
Let $M_{m,\Psi,\Phi}^{\left(E_1,E_2\right)}$ be a multiplier for $E_1$-Bessel sequence $\Phi=\lbrace\phi_k\rbrace_{k=1}^{\infty}\subseteq\mathcal{H}_1$ and $E_2$-Bessel sequence $\Psi=\lbrace\psi_k\rbrace_{k=1}^{\infty}\subseteq\mathcal{H}_2$. If $m\in\ell^\infty$, then $M_{m,\Psi,\Phi}^{\left(E_1,E_2\right)}$ is a well-defined bounded operator such that 
\begin{equation}
\left\Vert M_{m,\Psi,\Phi}^{\left(E_1,E_2\right)}\right\Vert\leq\sqrt{B}\sqrt{B^{'}}\left\Vert m\right\Vert_\infty.
\end{equation}
Moreover, $\sum_{k=1}^\infty m_k\left\langle f,\left(E\Phi\right)_k\right\rangle\left(E\Psi\right)_k$ converges unconditionally for all $f\in\mathcal{H}_1$.
\end{theorem}
\begin{proof}
By \eqref{eqmlt2} $M_{m,\Psi,\Phi}^{\left(E_1,E_2\right)}$ is well-defined because $\mathcal{M}_m$, $T_{E_2\Psi}$ and $T^\ast_{E_1\Phi}$ are well-defined by assumptions. By \cite[Lemma 5.4]{BAL1} $\mathcal{M}_m$ is bounded with $\left\Vert\mathcal{M}_m\right\Vert=\left\Vert m\right\Vert_\infty$. Therefore
\begin{equation*}
\left\Vert M_{m,\Psi,\Phi}^{\left(E_1,E_2\right)}\right\Vert=\left\Vert T_{E_2\Psi}\mathcal{M}_mT^\ast_{E_1\Phi}\right\Vert\leq\sqrt{B}\sqrt{B^{'}}\left\Vert m\right\Vert_\infty,
\end{equation*}
where $B$ and $B^{'}$ are $E_1$-Bessel and $E_2$-Bessel bounds for $\Phi$ and $\Psi$ respectively.
        
Finally, $\sum_{k=1}^\infty m_k\left\langle f,\left(E_1\Phi\right)_k\right\rangle\left(E_1\Psi\right)_k$ converges unconditionally for all $f\in\mathcal{H}_1$, since $\Psi$ is $E_2$-Bessel and $\left\lbrace m_k\left\langle f,\left(E_1\Phi\right)_k\right\rangle\right\rbrace_{k=1}^\infty$ belongs to $\ell^2(\mathbb{N})$ for all $f\in\mathcal{H}_1$ (see \cite{TAL}).                       
\end{proof}

 In the case $\mathcal{H}_1=\mathcal{H}_2$ and $E_1=E_2$, we use $M_{m,\Psi,\Phi}^{E}$ in notation.

\begin{theorem}
Suppose that $\Psi=\lbrace\psi_k\rbrace_{k=1}^{\infty}$ and $\Phi=\lbrace\phi_k\rbrace_{k=1}^{\infty}$ are dual frames for $\mathcal{H}$ and $E$ is any infinite diagonal complex matrix. Then the multipliers $M^{\overline{E},E^{-1}}_{1,\Psi,\Phi}$, $M^{\overline{E},E^{-1}}_{1,\Phi,\Psi}$, $M^{E^{-1},\overline{E}}_{1,\Psi,\Phi}$ and $M^{E^{-1},\overline{E}}_{1,\Phi,\Psi}$ are well defined. In fact, they are equal to $Id_{\mathcal{H}}$.
\end{theorem}
\begin{proof}
We just show that $M^{\overline{E},E^{-1}}_{1,\Phi,\Psi}=Id_{\mathcal{H}}$. The other cases will be proved analogously.

It is clear that $\Phi$ and $\Psi$ are $\overline{E}$-sequence and $E^{-1}$-sequence respectively. Furthermore, for given $f\in\mathcal{H}$, 
\begin{align*}
M^{\overline{E},E^{-1}}_{1,\Psi,\Phi}f&=\sum_{n=1}^{\infty}\left\langle f,\left(\overline{E}\left\lbrace\phi_k\right\rbrace_{k=1}^{\infty}\right)_n\right\rangle\left(E^{-1}\left\lbrace\psi_k\right\rbrace_{k=1}^{\infty}\right)_n\\&=\sum_{n=1}^{\infty}\left\langle f,\overline{E}_{n,n}\phi_n\right\rangle E^{-1}_{n,n}\psi_n\\&=\sum_{n=1}^{\infty}\left\langle f,\phi_n\right\rangle\psi_n=f.
\end{align*}
\end{proof}
\begin{lemma}\label{lem1}
Let $\lbrace\psi_k\rbrace_{k=1}^{\infty}$ be an $E$-frame for $\mathcal{H}$. Suppose that $\lbrace m_k\rbrace_{k=1}^{\infty}$ is a semi-normalized sequence of complex numbers. Then $\left\lbrace\left\vert m_k\right\vert\psi_k\right\rbrace_{k=1}^{\infty}$ is an $E$-frame for $\mathcal{H}$.
\end{lemma}
\begin{proof}
First note that $0<m=\inf_k\vert m_k\vert\leq\sup_k\vert m_k\vert=M<\infty$. Hence, for given $f\in\mathcal{H}$,
\begin{align*}
\sum_{n=1}^{\infty}\left\vert\left\langle f,\left(E\left\lbrace\left\vert m_k\right\vert\psi_k\right\rbrace_{k=1}^{\infty}\right)_n\right\rangle\right\vert^2&=\sum_{n=1}^{\infty}\left\vert\left\langle f,\sum_{k=1}^{\infty}E_{n,k}\left\vert m_k\right\vert\psi_k\right\rangle\right\vert^2\\&\leq M^2\sum_{n=1}^{\infty}\left\vert\left\langle f,\left(E\left\lbrace\psi_k\right\rbrace_{k=1}^{\infty}\right)_n\right\rangle\right\vert^2\leq M^2B\left\Vert f\right\Vert^2,
\end{align*}
where $B$ is the upper $E$-frame bound of $\left\lbrace\psi_k\right\rbrace_{k=1}^{\infty}$. Similarly, an straightforward argument shows that 
\begin{align*}
\sum_{n=1}^{\infty}\left\vert\left\langle f,\left(E\left\lbrace\left\vert m_k\right\vert\psi_k\right\rbrace_{k=1}^{\infty}\right)_n\right\rangle\right\vert^2\geq m^2A\Vert f\Vert^2, 
\end{align*}
where $A$ is the lower $E$-frame bound of $\left\lbrace\psi_k\right\rbrace_{k=1}^{\infty}$.
\end{proof}
\begin{remark}
Consider an infinite matrix $E$, a scalar sequence $m$ and a sequence $\Psi$ on $\mathcal{H}$. We want to examine the conditions under which
\begin{equation*}
m_n\left(E\left\lbrace\psi_k\right\rbrace_{k=1}^{\infty}\right)_n=\left(E\left\lbrace m_k\psi_k\right\rbrace_{k=1}^{\infty}\right)_n,
\end{equation*}
for all $n\in\mathbb{N}$. This means that
\begin{align*}
m_n\sum_{k=1}^\infty E_{n,k}\psi_k=\sum_{k=1}^\infty E_{n,k}m_k\psi_k. 
\end{align*}
But this is only true when 
\begin{equation}\label{eqmE}
m_nE_{n,k}=E_{n,k}m_k,
\end{equation}
for all $n,k\in\mathbb{N}$. As two particular examples, if $E$ is a diagonal matrix or $m$ is a constant sequence, \eqref{eqmE} will hold.
\end{remark}
The following theorem, shows that $M_{m,\Psi,\Psi}^{E}$ is bounded, when $\Psi$ is an $E$-frame.
\begin{theorem}
Consider an infinite matrix $E$ and a positive (resp. negative) and semi-normalized sequence $m=\{m_k\}$ which satisfy \eqref{eqmE}. Suppose that $\Psi=\left\lbrace\psi_k\right\rbrace_{k=1}^{\infty}$ is an $E$-frame on $\mathcal{H}$. Then $M_{m,\Psi,\Psi}^{E}=S_{E\lbrace\sqrt{m_k}\psi_k\rbrace}$ (resp. $M_{m,\Psi,\Psi}^{E}=-S_{E\lbrace\sqrt{\left\vert m_k\right\vert}\psi_k\rbrace}$) for the $E$-frame $\lbrace\sqrt{m_k}\psi_k\rbrace_{k=1}^{\infty}$ (resp. $\lbrace\sqrt{\left\vert m_k\right\vert}\psi_k\rbrace_{k=1}^{\infty}$) and is therefore invertible on $\mathcal{H}$.
\end{theorem}
\begin{proof}
As we show in Lemma \ref{lem1}, $\lbrace\sqrt{m_k}\psi_k\rbrace_{k=1}^{\infty}$ is an $E$-frame. \eqref{eqmE} implies that $E_{n,k}(m_n-m_k)=0$. So either $E_{n,k}=0$ or $m_n=m_k$ for all $n,k\in\mathbb{N}$. Equivalently, either $E_{n,k}=0$ or $\sqrt{m_n}=\sqrt{m_k}$ (resp. $\sqrt{\vert m_n\vert}=\sqrt{\vert m_k\vert}$) for positive (resp. negative) sequence $m$ and $\sqrt{m_n}E_{n,k}=E_{n,k}\sqrt{m_k}$ (resp. $\sqrt{\vert m_n\vert}E_{n,k}=E_{n,k}\sqrt{\vert m_k\vert}$) for all $n,k\in\mathbb{N}$.
Firstly, we prove the case where $m$ is positive.
For given $f\in\mathcal{H}$,
\begin{align*}
M_{m,\Psi,\Psi}^{E}f&=\sum_{n=1}^{\infty}m_n\left\langle f,\left(E\left\lbrace\psi_k\right\rbrace_{k=1}^{\infty}\right)_n\right\rangle\left(E\left\lbrace\psi_k\right\rbrace_{k=1}^{\infty}\right)_n\\&=\sum_{n=1}^{\infty}\left\langle f,\sqrt{m_n}\sum_{k=1}^\infty E_{n,k}\psi_k\right\rangle\sqrt{m_n}\sum_{k=1}^\infty E_{n,k}\psi_k\\&=\sum_{n=1}^{\infty}\left\langle f,\sum_{k=1}^\infty E_{n,k}\sqrt{m_k}\psi_k\right\rangle\sum_{k=1}^\infty E_{n,k}\sqrt{m_k}\psi_k\\&=\sum_{n=1}^{\infty}\left\langle f,\left(E\left\lbrace\sqrt{m_k}\psi_k\right\rbrace_{k=1}^{\infty}\right)_n\right\rangle\left(E\left\lbrace\sqrt{m_k}\psi_k\right\rbrace_{k=1}^{\infty}\right)_n=S_{E
\lbrace\sqrt{m_k}\psi_k\rbrace}f.
\end{align*}
If $m$ is a negative sequence, then $m_n=-\left(\sqrt{\left\vert m_n\right\vert}\right)^2$ for all $n$. Now, with a similar reasoning as the above argument, it is proven that $M_{m,\Psi,\Psi}^{E}f=-S_{E\lbrace\sqrt{\left\vert m_k\right\vert}\psi_k\rbrace}$.
\end{proof}
\begin{example}
Suppose that $\Psi=\left\lbrace e_1-e_2+\frac{1}{3}e_3, e_2-\frac{1}{3}e_3, \frac{1}{3}e_3, \frac{1}{4}e_4, \ldots\right\rbrace$ and 
\begin{equation*}
E=\begin{pmatrix}
1&1&0&0&\ldots\\
0&1&1&0&\ldots\\
0&0&3&0&\ldots\\
0&0&0&4&\ldots\\
\vdots&\vdots&\vdots&\vdots&\ddots
\end{pmatrix},
\end{equation*}which is invertible and the inverse is given by
\begin{equation*}
E^{-1}=\begin{pmatrix}
1&-1&\frac{1}{3}&0&\ldots\\
0&1&-\frac{1}{3}&0&\ldots\\
0&0&\frac{1}{3}&0&\ldots\\
0&0&0&\frac{1}{4}&\ldots\\
\vdots&\vdots&\vdots&\vdots&\ddots
\end{pmatrix}.
\end{equation*}
Then $E\Psi=\left\lbrace e_1, e_2, e_3, \ldots\right\rbrace.$ As $E_{1,2}$ and $E_{2,3}$ are non zero, set $m=\left\lbrace c, c, c, 1, \frac{1}{2}, \frac{1}{3}, \ldots\right\rbrace$ where $c>0$. Then $$m\left(E\Psi\right)=\left\lbrace ce_1, ce_2, ce_3, e_4, \frac{1}{2}e_5, \frac{1}{3}e_6, \ldots\right\rbrace.$$ On the other hand, $m\Psi=\left\lbrace ce_1-ce_2+\frac{c}{3}e_3, ce_2-\frac{c}{3}e_3, \frac{c}{3}e_3, \frac{1}{4}e_4, \frac{1}{10}e_5, \frac{1}{18}e_6, \ldots\right\rbrace$ and so we have $$E\left(m\Psi\right)=\left\lbrace ce_1, ce_2, ce_3, e_4, \frac{1}{2}e_5, \frac{1}{3}e_6, \ldots\right\rbrace.$$Hence $m(E\Psi)=E(m\Psi)$. Thus for given $f\in\mathcal{H}$,
\begin{align*}
M_{m,\Psi,\Psi}^{E}f&=\sum_{n=1}^\infty m_n\left\langle f,\left(E\left\lbrace\psi_k\right\rbrace_{k=1}^\infty\right)_n\right\rangle\left(E\left\lbrace\psi_k\right\rbrace_{k=1}^\infty\right)_n\\&= c\left\langle f,e_1\right\rangle e_1+c\left\langle f,e_2\right\rangle e_2+c\left\langle f,e_3\right\rangle e_3+\left\langle f,e_4\right\rangle e_4+\frac{1}{2}\left\langle f,e_5\right\rangle e_5+\frac{1}{3}\left\langle f,e_6\right\rangle e_6+\ldots.
\end{align*}
Since $m$ is semi-normalized, $\left\lbrace\sqrt{m_k}\psi_k\right\rbrace_{k=1}^\infty$ is an $E$-frame and $S_{E
\lbrace\sqrt{m_k}\psi_k\rbrace}$ is well-defined. Indeed,
$$\sqrt{m}\Psi=\left\lbrace \sqrt{c}e_1-\sqrt{c}e_2+\frac{\sqrt{c}}{3}e_3, \sqrt{c}e_2-\frac{\sqrt{c}}{3}e_3, \frac{\sqrt{c}}{3}e_3, \frac{1}{4}e_4, \frac{\sqrt{2}}{10}e_5, \frac{\sqrt{3}}{18}e_6, \ldots\right\rbrace$$
and $$E\left(\sqrt{m}\Psi\right)=\left\lbrace \sqrt{c}e_1, \sqrt{c}e_2, \sqrt{c}e_3, e_4, \frac{\sqrt{2}}{2}e_5, \frac{\sqrt{3}}{3}e_6, \ldots\right\rbrace.$$ Therefore
\begin{align*}
S_{E\lbrace\sqrt{m_k}\psi_k\rbrace}f&=\sum_{n=1}^\infty\left\langle f,\left(E\left\lbrace\sqrt{m_k}\psi_k\right\rbrace_{k=1}^\infty\right)_n\right\rangle\left(E\left\lbrace\sqrt{m_k}\psi_k\right\rbrace_{k=1}^\infty\right)_n\\&=\left\langle f, \sqrt{c}e_1\right\rangle\sqrt{c}e_1+\left\langle f, \sqrt{c}e_2\right\rangle\sqrt{c}e_2+\left\langle f, \sqrt{c}e_3\right\rangle\sqrt{c}e_3+\left\langle f,e_4\right\rangle e_4\\&+\left\langle f, \frac{\sqrt{2}}{2}e_5\right\rangle\frac{\sqrt{2}}{2}e_5+\left\langle f, \frac{\sqrt{3}}{3}e_6\right\rangle\frac{\sqrt{3}}{3}e_6+\ldots\\&=c\left\langle f,e_1\right\rangle e_1+c\left\langle f,e_2\right\rangle e_2+c\left\langle f,e_3\right\rangle e_3+\left\langle f,e_4\right\rangle e_4+\frac{1}{2}\left\langle f,e_5\right\rangle e_5+\frac{1}{3}\left\langle f,e_6\right\rangle e_6+\ldots,
\end{align*}
that is $M_{m,\Psi,\Psi}^{E}=S_{E\lbrace\sqrt{m_k}\psi_k\rbrace}$.
\end{example}
\begin{remark}
In \cite[Proposition 7.7]{BAL1} it is proven that if $\Psi=\lbrace \psi_k\rbrace_{k=1}^{\infty}$ and $\Phi=\lbrace \phi_k\rbrace_{k=1}^{\infty}$ are Riesz bases and $m=\lbrace m_n\rbrace_{n=1}^{\infty}$ is semi-normalized, then $M_{m,\Psi,\Phi}$ is invertible and $M^{-1}_{m,\Psi,\Phi}=M_{\frac{1}{m},\tilde{\Phi},\tilde{\Psi}}$, where $\tilde{\Phi}$ and $\tilde{\Psi}$ are the unique biorthogonal sequences of $\Phi$ and $\Psi$, respectively. The $E$-transform of an $E$-Riesz basis is a Riesz basis. Indeed,
\begin{align*}
\left(E\Psi\right)_n=\left(E\left\lbrace U\left(E^{-1}\left\lbrace e_j\right\rbrace_{j=1}^{\infty}\right)_k\right\rbrace_{k=1}^{\infty}\right)_n=\left(UE\left\lbrace\left(E^{-1}\left\lbrace e_j\right\rbrace_{j=1}^{\infty}\right)_k\right\rbrace_{k=1}^{\infty}\right)_n=Ue_n,
\end{align*} where $U$ is a bounded bijection on $\mathcal{H}$.
Thus we can extend the result mentioned above to $E$-Riesz beses.
\end{remark}
In the following, we can ask, when an $E$-multiplier is invertible, or when it is the inverse of another multiplier?
\begin{theorem}\label{th1}
Suppose that $E$ is an infinite matrix, $m=\{m_k\}$ is a positive (resp. negative) and semi-normalized sequence and \eqref{eqmE} holds. Let $\Phi=\left\lbrace\phi_k\right\rbrace_{k=1}^{\infty}$ be an $E$-frame and $\Psi=\left\lbrace\psi_k\right\rbrace_{k=1}^{\infty}$ be a sequence on $\mathcal{H}$ such that $\psi_k=U\phi_k$ for all $k\in\mathbb{N}$, where $U:\mathcal{H}\longrightarrow\mathcal{H}$ is a bounded bijection. Then $\Psi$ is an $E$-frame and $M_{m,\Phi,\Psi}^{E},M_{m,\Psi,\Phi}^{E}$ are invertible and
\begin{equation}\label{eqMS1}
\left(M_{m,\Phi,\Psi}^{E}\right)^{-1}=
\begin{cases}
\left(U^{-1}\right)^{\ast}S_{E\lbrace\sqrt{m_k}\psi_k\rbrace}^{-1}&m_k>0,\forall k,\\-\left(U^{-1}\right)^{\ast}S_{E\lbrace\sqrt{\left\vert m_k\right\vert}\psi_k\rbrace}^{-1}&m_k<0,\forall k.
\end{cases}
\end{equation}
\begin{equation}\label{eqMS2}
\left(M_{m,\Psi,\Phi}^{E}\right)^{-1}=
\begin{cases}
S_{E\lbrace\sqrt{m_k}\phi_k\rbrace}^{-1}U^{-1}&m_k>0,\forall k,\\-S_{E\lbrace\sqrt{\left\vert m_k\right\vert}\phi_k\rbrace}^{-1}U^{-1}&m_k<0,\forall k.
\end{cases}
\end{equation}
\end{theorem}
\begin{proof}
Firstly, we prove that $\Psi$ is an $E$-frame. For given $f\in\mathcal{H}$,
\begin{align*}
\sum_{n=1}^\infty\left\vert\left\langle f,\left(E\left\lbrace\psi_k\right\rbrace_{k=1}^\infty\right)_n\right\rangle\right\vert^2&=\sum_{n=1}^\infty\left\vert\left\langle f,\left(E\left\lbrace U\phi_k\right\rbrace_{k=1}^\infty\right)_n\right\rangle\right\vert^2\\&=\sum_{n=1}^\infty\left\vert\left\langle f,U\left (E\left\lbrace\phi_k\right\rbrace_{k=1}^\infty\right)_n\right\rangle\right\vert^2\leq B\Vert U\Vert^2\Vert f\Vert^2,
\end{align*}
where $B$ is an upper bound of $\Phi$. On the other hand,
\begin{align*}
\left\Vert f\right\Vert^2=\left\Vert\left(U^{-1}\right)^\ast U^\ast f\right\Vert^2\leq\left\Vert U^{-1}\right\Vert^2\left\Vert U^\ast f\right\Vert^2&\leq\frac{\left\Vert U^{-1}\right\Vert^2}{A}\sum_{n=1}^\infty\left\vert\left\langle U^\ast f,\left(E\left\lbrace\phi_k\right\rbrace_{k=1}^\infty\right)_n\right\rangle\right\vert^2\\&=\frac{\left\Vert U^{-1}\right\Vert^2}{A}\sum_{n=1}^\infty\left\vert\left\langle f,U\left(E\left\lbrace\phi_k\right\rbrace_{k=1}^\infty\right)_n\right\rangle\right\vert^2\\&=\frac{\left\Vert U^{-1}\right\Vert^2}{A}\sum_{n=1}^\infty\left\vert\left\langle f,\left(E\left\lbrace\psi_k\right\rbrace_{k=1}^\infty\right)_n\right\rangle\right\vert^2,
\end{align*}
where $A$ is a lower bound of $\Phi$. Hence $\Psi$ is an $E$-frame. 

Note that 
\begin{align*}
M_{m,\Phi,\Psi}^{E}f&=\sum_{n=1}^\infty m_n\left\langle f,\left(E\left\lbrace\psi_k\right\rbrace_{k=1}^\infty\right)_n\right\rangle\left(E\left\lbrace\phi_k\right\rbrace_{k=1}^\infty\right)_n\\&\sum_{n=1}^\infty m_n\left\langle U^\ast f,\left(E\left\lbrace\phi_k\right\rbrace_{k=1}^\infty\right)_n\right\rangle\left(E\left\lbrace\phi_k\right\rbrace_{k=1}^\infty\right)_n=M_{m,\Phi,\Phi}^{E}U^\ast f.
\end{align*}
By Theorem \ref{th1}, $M_{m,\Phi,\Phi}^{E}$ is invertible and 
\begin{equation*}
\left(M_{m,\Phi,\Phi}^{E}\right)^{-1}=
\begin{cases}
S_{E\lbrace\sqrt{m_k}\phi_k\rbrace}^{-1}&m_k>0,\forall k,\\-S_{E\lbrace\sqrt{\left\vert m_k\right\vert}\phi_k\rbrace}^{-1}&m_k<0,\forall k.
\end{cases}
\end{equation*}
Thus $M_{m,\Phi,\Psi}^{E}=M_{m,\Phi,\Phi}^{E}U^\ast$ is invertible and \eqref{eqMS1} is obtained. By similar argument, \eqref{eqMS2} also holds. 
\end{proof}
In the following results, we investigate some properties of $E$-multipliers.
\begin{theorem}
Let $H$ be a separable Hilbert space and $E=\left( E_{n,k}\right)_{n,k\geq1}$ an infinite complex diagonal matrix. Suppose that $\Psi=\lbrace \psi_k\rbrace_{k=1}^{\infty}$ and $\Phi=\lbrace \phi_k\rbrace_{k=1}^{\infty}$ are $E$-Bessel sequences for $\mathcal{H}$ and $m=\lbrace m_n\rbrace_{n=1}^{\infty}\in\ell^{\infty}$. If the sequence $\lbrace\lambda_n\rbrace_{n=1}^{\infty}=\lbrace E_{n,n}\rbrace_{n=1}^{\infty}$ belongs to $\ell^{\infty}$, then the multiplier operator for sequences $\lbrace\psi_n\rbrace_{n=1}^{\infty}$ and $\lbrace\phi_n\rbrace_{n=1}^{\infty}$ with symbol $m^{\prime}=\lbrace\vert\lambda_n\vert^2m_n\rbrace_{n=1}^{\infty}$ is well defined.
\end{theorem}
\begin{proof} 
The multiplier $M_{m, \Psi, \Phi}^{E}$ is well defined. Now for given $f\in\mathcal{H}$,
\begin{align*}
M_{m, \Psi, \Phi}^{E}f&=\sum_{n=1}^{\infty}m_n\left\langle f,\left(E\left\lbrace\phi_k\right\rbrace_{k=1}^{\infty}\right)_n\right\rangle\left(E\left\lbrace\psi_k\right\rbrace_{k=1}^{\infty}\right)_n\\&=\sum_{n=1}^{\infty}m_n\left\langle f,\lambda_n\phi_n\right\rangle\lambda_n\psi_n\\&=\sum_{n=1}^{\infty}\left\vert\lambda_n\right\vert^{2}m_n\left\langle f,\phi_n\right\rangle \psi_n=M_{\left(\left\vert\lambda_n\right\vert^{2}m_n\right),\Psi,\Phi}f.
\end{align*}
\end{proof}
\begin{theorem}
Consider two dual $E$-frames $\Psi=\lbrace\psi_k\rbrace_{k=1}^{\infty}$ and $\Phi=\lbrace\phi_k\rbrace_{k=1}^{\infty}$ in $\mathcal{H}$ where $E$ is an infinite complex diagonal matrix. Then $M_{m,\Psi,\Phi}=M_{m,\Phi,\Psi}=Id_{\mathcal{H}}$ with $m=\left\lbrace\left\vert E_{n,n}\right\vert^2\right\rbrace_{n=1}^{\infty}$.
\end{theorem}
\begin{proof}
For given $f\in\mathcal{H}$,
\begin{align*}
M_{m,\Psi,\Phi}f&=\sum_{n=1}^{\infty}\left\vert E_{n,n}\right\vert^2\left\langle f,\phi_n\right\rangle\psi_n\\&=\sum_{n=1}^{\infty}\left\langle f,E_{n,n}\phi_n\right\rangle E_{n,n}\psi_n\\&
=\sum_{n=1}^{\infty}\left\langle f,\left(E\left\lbrace\phi_k\right\rbrace_{k=1}^{\infty}\right)_n\right\rangle\left(E\left\lbrace\psi_k\right\rbrace_{k=1}^{\infty}\right)_n=f.
\end{align*}
Similarly, one can easily shows that $M_{m,\Phi,\Psi}=Id_{\mathcal{H}}$.
\end{proof}
We end this section with an example.
\begin{example}
Let $\mathcal{H}$ be a separable Hilbert space with an orthonormal basis $\lbrace e_k\rbrace_{k=1}^{\infty}$. Then $\Psi=\lbrace \psi_k\rbrace_{k=1}^{\infty}=\lbrace ke_k\rbrace_{k=1}^{\infty}$ and $\Phi=\lbrace \phi_k\rbrace_{k=1}^{\infty}=\lbrace k^2e_k\rbrace_{k=1}^{\infty}$
are non-Bessel sequences. Now consider the invertible infinite matrix $E$ as follow
\begin{equation*}
E_{n,j}=\begin{cases}
\frac{1}{n^2}&n=j,\\0&n\neq j.
\end{cases}
\end{equation*}
The matrix form of $E$ is 
\begin{equation*}
E=\begin{pmatrix}
1&0&0&\cdots\\
0&\frac{1}{4}&0&\cdots\\
0&0&\frac{1}{9}&\cdots\\
\vdots&\vdots&\vdots&\ddots
\end{pmatrix}.
\end{equation*} 
An easy argument shows that $\lbrace \psi_k\rbrace_{k=1}^{\infty}$ and $\lbrace \phi_k\rbrace_{k=1}^{\infty}$ are $E$-Bessel sequences. Thus $M^{E}_{m,\Psi,\Phi}$ is well defined for all $m\in\ell^{\infty}$. Suppose that $m=1$ and $f\in\mathcal{H}$. Then 
\begin{align*}
M_{m,\Psi,\Phi}^{E}f&=\sum_{n=1}^{\infty}1\left\langle f,\left(E\left\lbrace\phi_k\right\rbrace_{k=1}^{\infty}\right)_n\right\rangle\left(E\left\lbrace\psi_k\right\rbrace_{k=1}^{\infty}\right)_n\\&=\sum_{n=1}^{\infty}\left\langle f,\frac{1}{n^2}\phi_n\right\rangle\frac{1}{n^2}\psi_n\\&=\sum_{n=1}^{\infty}\frac{1}{n^4}\left\langle f,\phi_n\right\rangle \psi_n=M_{\left(\frac{1}{n^4}\right),\Psi,\Phi}f\\&=\sum_{n=1}^{\infty}\frac{1}{n}\left\langle f,e_n\right\rangle e_n=M_{\left(\frac{1}{n}\right),\left(e_n\right),\left(e_n\right)}f.
\end{align*}
\end{example}


\begin{thebibliography}{0}
\bibitem{BAL1}  P. Balazs, \emph{Basic definition and properties of Bessel multipliers},  J. Math. Anal. Appl. 325 (2007) 571–585.

\bibitem{BAL2}  P. Balazs, D.T.Stoeva, \emph{Detailed characterization of unconditional convergence and invertibility of multipliers}, Sampling Theory in Signal and Image Processing, \textbf{12} (2013), 87–-125.

\bibitem{BAL3} P. Balazs, B. Laback, G. Eckel and W.A. Deutsch,  \emph{Time-Frequency Sparsity by
Removing Perceptually Irrelevant Components Using a Simple Model of Simultaneous Masking}, IEEE Transactions on Audio, Speech and Language Processing \textbf{18}(1) (2010), 34–-49.

\bibitem{CA} P.G. Casazza, \emph{The art of frame theory}, Taiwan. J. Math., \textbf{4}(2) (2000), 129-–201.

\bibitem{ch} O. Christensen, \emph{An introduction to frames and Riesz bases}, Birkhauser, Boston, 2016.

\bibitem{Conway} J. B. Conway, \emph{A Course in Functional Analysis}, 2nd Edition, Springer-Verlag, 1990.

\bibitem{DG} I. Daubechies, A. Grassman and Y. Meyer, \emph{Painless nonothogonal expanisions}, J. Math. Phys.,  
\textbf{27} (1986), 1271--1283.

\bibitem{DS} R.J. Duffin and A.C. Schaeffer, \emph{A class of nonharmonic Fourier series}, Trans. Amer.
Math. Soc., \textbf{72} (1952), 341--366.

\bibitem{HAN} D. Han and D. Larson, \emph{Frame, bases and group representations}, Memoir. Amer. Math. Soc., \textbf{147} (2000), 1–-94.

\bibitem{MAJ1} P. Majdak, P. Balazs and B. Laback,  \emph{Multiple Exponential Sweep Method for Fast
Measurement of Head Related Transfer Functions}, J. Audio Engineering Society, \textbf{55}(7/8) (2007), 623–-637.

\bibitem{MAJ2} P. Majdak, P. Balazs, W.Kreuzer and M. Dorfler, \emph{A Time-Frequency Method
for Increasing the Signal-To-Noise Ratio in System Identification with Exponential
Sweeps}, Proceedings of the 36th International Conference on Acoustics, Speech and Signal Processing, ICASSP 2011, Prag, 2011.

\bibitem{TAL} G. Talebi, M. A. Dehghan, \emph{On $E$-frames in separable Hilbert spaces},
Banach J. Math. Anal. 9(3) (2015), 43--74.



\end{thebibliography}
\end{document}